\pdfoutput=1 

\documentclass{svmult}

\usepackage{graphicx}
\usepackage{multicol}        		
\usepackage[bottom]{footmisc}		
\usepackage{newtxtext}
\usepackage[varvw]{newtxmath}       

\usepackage{enumitem}
\usepackage{mathtools}

\usepackage{adjustbox}
\usepackage{caption}
\usepackage{subcaption}
\usepackage{tikz-cd} 
\usetikzlibrary{arrows.meta,calc,positioning,fit,backgrounds,shapes.geometric}
\usepackage{adjustbox}

\newtheorem{assumption}{Assumption}

\usepackage{xcolor}
\usepackage[colorlinks]{hyperref}
\usepackage{bookmark}

\colorlet{siaminlinkcolor}{green!50!black}
\colorlet{siamexlinkcolor}{red!50!black}
\colorlet{siamreviewcolor}{black!50}
\hypersetup{allcolors=siaminlinkcolor,urlcolor=siamexlinkcolor}

\ifpdf
\hypersetup{
	pdftitle={A short tour of operator learning theory: Convergence rates, statistical limits, and open questions},
	pdfauthor={S. Brugiapaglia, N. R. Franco, and N. H. Nelsen}
}
\fi

\makeatletter

\def\toclevel@title{1000}
\def\toclevel@titlech{1000}
\def\toclevel@author{1000}
\def\toclevel@authorch{1000}

\providecommand*\toclevel@bibtoc{1000} 
\let\l@bibtoc\l@section               

\AtBeginDocument{%
	\@ifundefined{bibsection}{}{%
		\def\bibsection{%
			\section*{\refname}
			\addcontentsline{toc}{bibtoc}{\refname}
			\csname biblst@rthook\endcsname\par
		}%
	}%
}

\makeatother

\begin{document}

\title*{A short tour of operator learning theory:\\
Convergence rates, statistical limits,\\
and open questions}

\titlerunning{A short tour of operator learning theory}
\author{Simone Brugiapaglia\orcidID{0000-0003-1927-8232},\\Nicola Rares Franco\orcidID{0000-0002-4980-5606},\\Nicholas H. Nelsen\orcidID{0000-0002-8328-1199}}
\institute{Simone Brugiapaglia \at Concordia University, \email{simone.brugiapaglia@concordia.ca}
\and Nicola Rares Franco \at Politecnico di Milano, \email{nicolarares.franco@polimi.it}
\and Nicholas H. Nelsen \at Cornell University and The University of Texas at Austin, \email{nnelsen@oden.utexas.edu}}

\maketitle

\abstract*{\indent
This paper surveys recent developments at the intersection of operator learning, statistical learning theory, and approximation theory. First, it reviews error bounds for empirical risk minimization with a focus on holomorphic operators and neural network approximations. Next, it illustrates fundamental performance limits in terms of sample size by adopting a minimax perspective and considering various notions of regularity beyond holomorphy. The paper ends with a discussion on the interplay between these two perspectives and related open questions.
}

\abstract{\indent
This paper surveys recent developments at the intersection of operator learning, statistical learning theory, and approximation theory. First, it reviews error bounds for empirical risk minimization with a focus on holomorphic operators and neural network approximations. Next, it illustrates fundamental performance limits in terms of sample size by adopting a minimax perspective and considering various notions of regularity beyond holomorphy. The paper ends with a discussion on the interplay between these two perspectives and related open questions.
}

\section{Introduction}\label{sec:intro}
Since their first introduction, our theoretical understanding of deep neural networks and related architectures has increased remarkably. Already in the 1990s, Cybenko~\cite{cybenko1989approximation} and Hornik~\cite{hornik1991approximation} proved several foundational results, referred to as \emph{universal approximation} theorems, which established the density of neural networks in classical function spaces such as $C^k$ and $L^p$ spaces. Relatedly, Chen~and~Chen~\cite{chen1995universal} extended these arguments to the case of nonlinear operators between function spaces by introducing a neural network-based architecture later referred to as \emph{DeepONet} \cite{lu2021learning}.

These findings about existence \cite{lanthaler2022error} were later refined with results on model \emph{expressivity}. This established a connection between architectural complexity and model accuracy \cite{guhring2021approximation,marcati2023exponential,petersen2018optimal,schwab2023deep,siegel2023optimal}. See, e.g., work by Mhaskar~\cite{mhaskar1996neural} on smooth shallow networks and Yarotsky~\cite{yarotsky2017error} on deep ReLU architectures.

However, all such results leave out an important part of deep learning: the \emph{training} process. In practice, neural networks are typically trained over empirical data, and they can only learn from a limited amount of samples. Additionally, the training procedure usually involves the highly nontrivial minimization of a strongly nonconvex objective functional---the \emph{loss function}---making the problem even more challenging.

In this sense, a comprehensive picture of deep learning for scientific computing requires a proper interplay of approximation theory, statistical learning, and optimization theory. This paper focuses on the second pillar in the context of \emph{operator learning} \cite{boulle2024mathematical,kovachki2024operator,nelsen2024operator,nelsen2025operator,subedi2025operator}.
Specifically, the present work aims to summarize recent contributions in this field, discuss their implications, and highlight outstanding challenges.
These considerations include error bounds and sample complexity rates for learning special classes of target operators, such as holomorphic ones, but also more general results holding under minimal regularity assumptions, such as Lipschitz continuity; see Fig.~\ref{fig:flowchart} for a schematic overview. In the former case, the paper reviews and compares explicit constructions proposed in the literature, emphasizing the different analytical techniques employed in their derivation. In the latter case, instead, the paper presents results that hold in the \emph{minimax} sense and encompass both lower and upper rates. The minimax framework raises fundamental questions about the tightness of the holomorphic analysis and the necessity of properly characterizing the operator classes of interest, potentially beyond sole regularity requirements.

\runinhead{Outline.}
The paper is organized as follows. \S\ref{sec:erm} presents two sample complexity results for supervised training of neural operators using empirical risk minimization (ERM) \cite{gyorfi2002distribution}, whereas \S\ref{sec:minimax} comments on the optimality of these results from the minimax perspective \cite[Chp.~9]{krieg2026optimal}. \S\ref{sec:conclusion} summarizes the findings and collects open problems.

\begin{figure}[tb]
	\centering
	\captionsetup{skip=10pt}
	
	\definecolor{olblue}{HTML}{0E7490}
	\definecolor{olbluedark}{HTML}{083F5B}
	\definecolor{ollight}{HTML}{E8F5F8}
	\definecolor{olgray}{HTML}{F5F7FA}
	\definecolor{olorange}{HTML}{FFF3E0}
	\definecolor{darkgreen}{rgb}{0,0.6,0}

    \definecolor{lightblue}{rgb}{0.6,0.8,1}

	\begin{adjustbox}{max width=\textwidth}
	\begin{tikzpicture}[
		x=1cm,y=1cm,
		font=\sffamily\scriptsize,
		arrow/.style={-{Stealth[length=1.8mm]}, line width=0.85pt, draw=black!80},
		softarrow/.style={-{Stealth[length=1.8mm]}, line width=0.85pt, dashed, draw=black!55},
		linkarrow/.style={-{Stealth[length=1.8mm]}, line width=0.85pt, dotted, draw=olbluedark!85},
		main/.style={rectangle, rounded corners=5pt, draw=olbluedark, fill=olbluedark, text=white, align=center, text width=1.95cm, minimum height=0.9cm, inner sep=4pt, font=\sffamily\scriptsize\bfseries},
		class/.style={rectangle, rounded corners=4pt, draw=olbluedark!65, fill=ollight, align=center, text width=1.9cm, minimum height=0.55cm, inner sep=3pt},
		openclass/.style={rectangle, rounded corners=4pt, draw=orange!75!black, fill=olorange, align=center, text width=2.15cm, minimum height=0.55cm, inner sep=3pt},
		result/.style={rectangle, rounded corners=4pt, draw=black!70, fill=white, align=left, text width=3.05cm, minimum height=0.82cm, inner sep=4pt},
		wideopen/.style={rectangle, rounded corners=4pt, draw=orange!80!black, fill=olorange, align=left, text width=3.35cm, inner sep=4pt},
		group/.style={rectangle, draw=none, fill=white, inner sep=0pt},
		tag/.style={font=\sffamily\tiny\itshape, text=black!60, align=center, fill=white, fill opacity=0.9, text opacity=1, inner sep=1pt}
		]
				
		\node[main] (erm) at (0,4.35) {{\normalfont{\S\hypersetup{linkcolor=lightblue}\ref{sec:erm}}}: Empirical risk\\minimization};
		\node[class] (holoERM) at (2.65,4.80) {Holomorphic\\operators};
		\node[class] (model) at (5.20,4.80) {\textbf{Assumption}~\ref{assumption:setting}\\Encoder-decoder neural operators\\
        \tiny{e.g., FrameNet, DeepONet, ...}};
		\node[class] (roughERM) at (2.65,3.75) {Less regular\\operators};
		
		\node[result] (thm1) at (8.25,5.05) {
			\textbf{Theorem~\ref{thm:upper_ep}}~\cite{reinhardt2024statistical}\\
			\emph{Near Monte Carlo upper rate}\\
			Fully trainable; random noise
		};
		
		\node[result] (thm2) at (8.25,3.65) {
			\textbf{Theorem~\ref{thm:upper_cs}}~\cite{adcock2024optimal}\\
			\emph{Fast algebraic upper rates}\\
			Partially trainable;
            noise uniformly bounded by $\sigma \geq 0$
		};
		
		\node[wideopen] (openERM) at (12.45,4.35) {
			\textbf{Open questions}\vspace{0.075cm}\\
			Can the ERM Monte Carlo\\barrier be broken?\vspace{0.075cm}\\
            Fast rates for fully trainable\\networks?
		};
		
		\draw[arrow] (erm.east) -- (holoERM.west);
		\draw[arrow] (erm.east) -- (roughERM.west);
		\draw[arrow] (holoERM.east) -- (model.west);
		\draw[arrow] (model.east) -- (thm1.west);
		\draw[arrow] (model.east) -- (thm2.west);
		\draw[arrow] (thm1.east) -- (openERM.west);
		\draw[arrow] (thm2.east) -- (openERM.west);
		
		\node[main] (minimax) at (0,0.35) {{\normalfont\S\hypersetup{linkcolor=lightblue}\ref{sec:minimax}}: General learning methods\\(minimax theory)};
		\node[class] (widths) at (2.65,0.35) {Minimax lower and upper bounds};
		\node[class] (holoClass) at (5.20,1.75) {Holomorphic operators};
		\node[class] (fnoClass) at (5.20,0.35) {Operators well-approximated by FNOs};
		\node[class] (lipClass) at (5.20,-1.05) {Lipschitz and\\$C^k$ balls};
		
		\node[result] (thm4) at (8.25,1.75) {
			\textbf{Theorem~\ref{thm:minimax_holo}}~\cite{adcock2024optimal}\\
			\emph{Noise-free fast algebraic lower rates}
		};
		
		\node[result] (thm5) at (8.25,0.35) {
			\textbf{Theorem~\ref{thm:fno_optimal_rate}}~\cite{grohs2025theory,kovachki2024data}\\
			\emph{Two-sided bounds on best possible algebraic upper rate}
		};
		
		\node[result] (thm3) at (8.25,-1.05) {
			\textbf{Theorem~\ref{thm:lower_Ck}}~\cite{kovachki2024data}\\
			\emph{Algebraic minimax rates are\\impossible; polylog instead}
		};
		
		\node[result] (thm6) at (12.45,-1.05) {
			\textbf{Theorem~\ref{thm:lip_upper_lower}}~\cite{adcock2025towards,kovachki2024data}\\
			\emph{Polylog scaling of \emph{noisy} sampling $n$-width is unavoidable}
		};
		
		\node[wideopen] (openMM) at (12.45,1.10) {
			\textbf{Open questions}\\\vspace{0.075cm}
			Is ERM minimax optimal\\beyond holomorphy?\\\vspace{0.075cm}
			Sharp exponent in Theorem~\ref{thm:fno_optimal_rate}?\\\vspace{0.075cm}
            Noise level aware minimax rates for various classes?
		};
		
		\draw[arrow] (minimax.east) -- (widths.west);
		\draw[arrow] (widths.east) -- (holoClass.west);
		\draw[arrow] (widths.east) -- (fnoClass.west);
		\draw[arrow] (widths.east) -- (lipClass.west);
		\draw[arrow] (holoClass.east) -- (thm4.west);
		\draw[arrow] (fnoClass.east) -- (thm5.west);
		\draw[arrow] (lipClass.east) -- (thm3.west);
		\draw[arrow] (thm3.east) -- node[tag, above, yshift=2.5pt] {noisy\\analog} (thm6.west);
		\draw[arrow] (thm4.east) -- (openMM.west);
		\draw[arrow] (thm5.east) -- (openMM.west);
		\draw[arrow] (thm6.north) -- (openMM.south);

		\draw[softarrow] (roughERM.south) to[out=-90,in=90] node[tag, left, xshift=-2pt] {ERM upper\\bounds} (widths.north);
		\draw[softarrow] (thm2.south) to[out=-90,in=90] node[tag, right, xshift=1pt] {matches minimax optimal\\rate when $\sigma=0$} (thm4.north);
		\draw[softarrow, 	{Stealth[length=1.8mm]}-{Stealth[length=1.8mm]}] (openERM.south) to[out=-90,in=90] (openMM.north);
	\end{tikzpicture}
	\end{adjustbox}

	\caption{High-level overview of the empirical risk minimization (ERM) upper bounds reviewed in \S\ref{sec:erm}, minimax lower bounds from \S\ref{sec:minimax}, and open questions from \S\ref{sec:conclusion}.
	Solid arrows indicate logical flow within each perspective. Dashed arrows indicate
	bridges between upper and lower bounds.}
	\label{fig:flowchart}
\end{figure}

\section{Error bounds for empirical risk minimization}\label{sec:erm}
In this section, we overview two recent results on operator learning that focus on ERM via deep learning, namely \cite{reinhardt2024statistical} and \cite{adcock2024optimal}. Both contributions address the case of \emph{holomorphic} operators, thus requiring strong regularity assumptions. This choice is supported by several reasons. First, holomorphy emerges naturally in the study of solution operators of parameter-dependent partial differential equations (PDEs) \cite{cohen2015approximation, cohen2018shape}, where operator learning plays a central role. Second, the rich mathematical structure of holomorphic operators enables a more refined analysis of both analytical and statistical aspects, offering a strong foundation for future research. In this regard, it is noteworthy that the learning problem associated with holomorphic operators can be approached using quite different mathematical techniques; as we shall see, the two contributions summarized hereafter provide a clear illustration of such versatility.
The former, \cite{reinhardt2024statistical}, leverages the theory of \emph{empirical processes},
whereas the latter, \cite{adcock2024optimal}, is grounded in the theory of \emph{compressed sensing}. 

For ease of exposition and in an attempt to unify the setting and notation, we adopt a different presentation of the results than the one by the original authors. Clearly, this process inevitably results in a partial loss of generality. The interested reader is invited to consult the original works for additional details and more general formulations of the results showcased here.
Throughout \S\ref{sec:erm}, we consider the following assumptions.

\begin{assumption}
    \label{assumption:setting}
    Let $\mathscr{G}\colon \mathcal{U}\to \mathcal{V}$ be a (nonlinear) operator between two real separable Hilbert spaces of infinite dimension. Let $\varrho$ be a probability measure over $\mathcal{U}$. Let $\mathscr{E}_\infty\colon\mathcal{U}\to \ell^2(\mathbb{N})$ and $\mathscr{D}_\infty\colon\ell^2(\mathbb{N})\to\mathcal{V}$ both be boundedly invertible linear maps. For each $s\in\mathbb{N}$, denote by  $\mathscr{R}_s\colon\ell^2(\mathbb{N})\to\mathbb{R}^s$ the restriction operator mapping a sequence onto its first $s$ entries. Write $\mathscr{R}_s^*$ for its adjoint, i.e., the ``zero padding'' operator. 
\end{assumption}

The goal of operator learning is to approximate $\mathscr{G}$ by means of a suitable neural operator architecture trained on labeled---and potentially noisy---data. Among many possible classes of neural operators, we focus on architectures of the form
\begin{align*}
    \hat{\mathscr{G}}=\mathscr{D}_q\circ g \circ \mathscr{E}_d,\quad\text{where}
\end{align*}

\begin{enumerate}[label=(\emph{\roman*}),leftmargin=1.75\parindent,itemsep=0pt, labelsep=4pt]
    \item $\mathscr{E}_d\colon \mathcal{U}\to\mathbb{R}^{d}$ is
the \emph{encoder} map, given as $\mathscr{E}_d\coloneqq\mathscr{R}_d\circ \mathscr{E}_\infty$;
    \item $\mathscr{D}_q\colon \mathbb{R}^{q}\to\mathcal{V}$ is
the \emph{decoder} map, defined as $\mathscr{D}_q\coloneqq\mathscr{D}_\infty\circ \mathscr{R}_q^*$;
    \item $g\colon\mathbb{R}^{d}\to\mathbb{R}^{q}$ is a neural network model, such as a multilayer perceptron (MLP).
\end{enumerate}

Essentially, the encoder and the decoder serve the purpose of handling the infinite dimensions characterizing $\mathcal{U}$ and $\mathcal{V}$, whereas $g$ seeks to approximate $\mathscr{G}$ at the latent level.
This setup includes, for instance, architectures such as \emph{PCA-Net} \cite{bhattacharya2021model} and \emph{Frame-Net} \cite{reinhardt2024statistical}; furthermore, it shares interesting connections with other popular techniques such as DeepONets \cite{lu2021learning} and Fourier Neural Operators (FNOs) \cite{li2020fourier}.

In the setting considered here, the MLP $g$ is the only trainable component of the neural operator.
The truncation indices for the encoder and the decoder, $d$ and $q$, are typically chosen based on the amount of data at hand. They usually reflect regularity properties associated to the elements in the domain and codomain of $\mathscr{G}$. 

Given a finite collection of noisy samples $\{(u_i,v_i)\}_{i=1}^n\subset\mathcal{U}\times\mathcal{V}$, where
\begin{equation}
    \label{eq:labelled data}
u_i\stackrel{\mathrm{i.i.d.}}{\sim}\varrho\quad\text{and}\quad v_i\coloneqq\mathscr{G}(u_i)+e_i\quad\text{with}\quad e_i \in \mathcal{V},
\end{equation}
the objective is to construct an approximation $\hat{\mathscr{G}}$ such that the mean squared error
\begin{equation}
\label{eq:mse}
\|\mathscr{G}-\hat{\mathscr{G}}\|_{L^2_{\varrho}(\mathcal{U};\mathcal{V})}^2\coloneqq \int_{\mathcal{U}}\|\mathscr{G}(u)-\hat{\mathscr{G}}(u)\|_{\mathcal{V}}^2\,d\varrho(u)
\end{equation}
is minimized. 
Here, $e_i \in \mathcal{V}$ represents external noise. It can be treated either as deterministic (worst-case analysis, e.g., discretization error) or stochastic (e.g., statistical noise). In
ERM, the learning procedure consists in solving the optimization problem
\begin{equation}
\label{eq:training_pb}
\min_{g\in\mathscr{N}}\;\frac{1}{n}\sum_{i=1}^{n}\|v_i-(\mathscr{D}_q\circ g\circ \mathscr{E}_d)(u_i)\|_{\mathcal{V}}^2,
\end{equation}
where $\mathscr{N}$ is a suitable neural network class that depends on $d$ and $q$. Given a minimizer $\hat{g}_n \in \mathscr{N}$ of \eqref{eq:training_pb}, the associated operator approximation is defined as $\hat{\mathscr{G}}_n \coloneqq \mathscr{D}_q \circ \hat{g}_n \circ \mathscr{E}_d$.

\runinhead{Empirical process approach.}
In \cite{reinhardt2024statistical}, the authors derive statistical error bounds for operator learning of holomorphic operators by combining previously established results from approximation theory, mostly \cite{herrmann2024neural,schwab2023deep}, and empirical process theory \cite{geer2000empirical,vaart2023empirical}. The latter provides powerful tools for analyzing the statistical behavior of empirical risk minimizers and has been exploited extensively in the literature \cite{lanthaler2023error,liu2024deep}. Up to some simplifications, which serve to ease notation, we can summarize one of the main results of \cite{reinhardt2024statistical} as follows.
Given $r> 1/2$, $R>0$, and $t\geq0$, let
$$
\mathsf{H}_{r,R}\coloneqq\prod_{k\in\mathbb{N}}[-Rk^{-r}, Rk^{-r}]\subset \ell^2(\mathbb{N})\quad\text{and}\quad\|x\|_{2,t}\coloneqq\Biggl(\sum_{k=1}^{\infty}k^{2t}x_k^2\Biggr)^{1/2}
$$
denote a \emph{Hilbert-type hypercube} of power $r$ and a $t$-weighted $\ell^2$ norm, respectively. Given any $N\in\mathbb{N}$ and $c'>0$, we denote by $\mathscr{N}_{N,c'}$ the class of deep ReLU networks whose width does not exceed $c'N$ and whose depth is smaller than $c'\log(N)$. As a corollary of \cite[Thm.~3.15]{reinhardt2024statistical}, we have the following result. 

\begin{theorem}\label{thm:upper_ep}
    \label{th:jakob}
    Instate Assumption~\ref{assumption:setting}. Assume that for some $r>1$ and $R>0$, the operator $\mathscr{G}$ admits an holomorphic extension $\tilde{\mathscr{G}}$ over a complex open set $\mathcal{O}$ containing $\mathscr{E}^{-1}_{\infty}(\mathsf{H}_{r,R})$.  Further assume that for some $t>0$,
    $$
    \mathrm{supp}(\varrho) \subseteq \mathscr{E}^{-1}_{\infty}(\mathsf{H}_{r,R})\quad \text{and}\quad
    \sup_{u\in \mathcal{O}}\|\mathscr{D}_{\infty}^*\mathscr{G}(u)\|_{2,t}<\infty.
    $$
    Let $\kappa\coloneqq2\min\{r-1,t\}$ and fix $\tau>0$ arbitrarily small.
    Then there exist positive constants $c$ and $c'$ such that the following holds.
    For every sample size $n\in\mathbb{N}$ and every training set $\{(u_i,v_i)\}_{i=1}^{n}$ as in \eqref{eq:labelled data} with the $\{e_i\}_{i=1}^n$ i.i.d. subgaussian and independent of the $\{u_i\}_{i=1}^n$,
    the class  $\mathscr{N} = \mathscr{N}_{N,c'}$ obtained by setting $N=\lceil n^{1/(\kappa+1)}\rceil$ and $d = q = N(n)$ admits at least one minimizer  $\hat{g}_n\in\mathscr{N}$  of \eqref{eq:training_pb}. Additionally, for any such minimizer,
    the corresponding trained neural operator  $\hat{\mathscr{G}}_n=\mathscr{D}_{q(n)}\circ \hat{g}_n\circ \mathscr{E}_{d(n)}$ satisfies
    \begin{equation}
    \label{eq:jakob-bound}
    \Bigl(\mathbb{E}\|\mathscr{G}-\hat{\mathscr{G}}_n\|_{L^2_{\varrho}(\mathcal{U};\mathcal{V})}^2
    \Bigr)^{1/2}
    \le c n^{-\frac{1}{2}\left(\frac{1}{1+2/\kappa}\right)+\tau},
    \end{equation}
    where the expectation is taken with respect to the law of the $\{u_i\}_{i=1}^{n}$ and $\{e_i\}_{i=1}^n$.
\end{theorem}

The error bound in \eqref{eq:jakob-bound} exhibits an approximate Monte Carlo rate when $r\uparrow\infty$, $t\uparrow\infty$, and $\tau\downarrow0$, with the right-hand side approaching $n^{-1/2}$. However, the actual decay ultimately depends on the values of the three exponents $r$, $t$, and $\tau$, whose interpretation is given as follows.
The former, $r$ and $t$, are related to the regularity of the elements in the input and output spaces, respectively. For instance, when $\mathcal{U}$ and $\mathcal{V}$ are function spaces, the two can describe how rapidly the spectrum of the elements $u\in\text{supp}(\varrho)\subseteq\mathcal{U}$ and $v\in \mathscr{G}(\text{supp}(\varrho))\subseteq\mathcal{V}$ decay with respect to a chosen basis. 

Conversely, the $\tau$ exponent is merely an algebraic artifact introduced to suppress log factors in $n$. The latter, in fact, commonly appear when relying on statistical learning arguments that decompose the mean squared error \eqref{eq:mse} into an \emph{approximation error} and a \emph{statistical error} associated with the expressivity and the metric entropy of the model class, respectively. In \cite{reinhardt2024statistical}, the authors characterize the approximation error by carefully balancing truncation errors---arising from the encoder-decoder design---with previously known results on the exponential convergence of neural network approximation for holomorphic functions. The entropy of ReLU network classes, instead, is analyzed explicitly by leveraging properties of the ReLU activation function.

Compared to Thm.~\ref{th:jakob}, the analysis presented in \cite{reinhardt2024statistical} is more general in at least two directions. First, the authors also provide generalization bounds for sparse networks, i.e., for a specialized subclass $\mathscr{N}^{(s)}\subseteq\mathscr{N}$. Second, they consider scenarios in which the encoding and decoding procedures can arise from arbitrary \emph{frames}, whereas our construction is only compatible with representations derived from a Riesz basis.

\runinhead{Compressed sensing approach.} 
We now present the error analysis from \cite{adcock2024optimal}. We start by defining the holomorphic regularity considered therein. For  $b \in \ell^1(\mathbb{N})$ with $b_j \geq 0$ for each $j$, define the set
\begin{align*}
    \mathsf{BP}(b)\coloneqq \biggl\{\,\mathcal{P}(\rho) \,\,\Big|\,\,\rho\in\mathbb{R}^{\mathbb{N}} \text{ such that } \rho_j \geq 1 \text{ for every } j, \; \sum_{j=1}^\infty \biggl(\frac{\rho_j + \rho_j^{-1}}{2}-1\biggr) b_j \leq 1\,\biggr\}
\end{align*}
of \emph{filled-in Bernstein polyellipses} $\mathcal{P}(\rho) \coloneqq \mathcal{B}(\rho_1) \times \mathcal{B}(\rho_2) \times \cdots \subseteq \mathbb{C}^{\mathbb{N}}$, where $\mathcal{B}(\rho_j) \coloneqq \{(z + z^{-1})/2 \,|\, z \in \mathbb{C},\, 1\leq |z|\leq \rho_j\}$. Then, define the region $\mathcal{R}(b) \subseteq \mathbb{C}^{\mathbb{N}}$ by the union $\mathcal{R}(b) \coloneqq \bigcup_{\mathcal{P}\in \mathsf{BP}(b)}\mathcal{P}$ of all filled-in Bernstein polyellipses. Note that $\mathsf{H}_{0,1} = [-1,1]^{\mathbb{N}} \subseteq \mathcal{R}(b)$.
This leads to the function class 
\begin{equation}
\label{eq:def_H(b)}
\mathcal{H}(b) \coloneqq \bigl\{f\colon \mathsf{H}_{0,1} \to \mathcal{V} \,\big|\, f = \tilde{f}|_{\mathsf{H}_{0,1}}, \, \tilde{f} \text{ holomorphic in } \mathcal{R}(b)\text{, }\|\tilde{f}\|_{L^{\infty}(\mathcal{R}(b); \mathcal{V})} \leq 1\bigr\}.
\end{equation}
As aforementioned, solution operators to parametric PDEs such as the parametric diffusion equation with affinely-parametrized diffusion coefficient belong to the class $\mathcal{H}(b)$; see, e.g., \cite[Prop.~4.9]{adcock2022sparse}. Below, we consider sequences $b \in \ell^p_{\mathsf{M}}(\mathbb{N})$; this means that $\tilde{b} \in \ell^p(\mathbb{N})$, where $\tilde{b}_j = \sup_{i \geq j} |b_i|$ for each $j$.

We now present a corollary of \cite[Thm.~3.1]{adcock2024optimal}. Specifically, here we assume the map $\iota$ from \cite{adcock2024optimal} coincides with the boundedly invertible linear operator $\mathscr{E}_\infty$ of Assumption~\ref{assumption:setting}. 

\begin{theorem}\label{thm:upper_cs}
Let $n\geq 3$, $\delta >0$ arbitrarily small, and set $\tilde{n} \coloneqq n/\log^4(n)$. Then there exists a class $\mathscr{N}$ of tanh neural networks depending solely on $n$, with 
\[\mathrm{width}(\mathscr{N}) \le c_0\, \tilde{n}^{1+\delta} \quad \text{and} \quad \mathrm{depth}(\mathscr{N}) \le c_0 \log(\tilde{n})
\]
for $c_0>0$ a suitable universal constant,
such that the following holds. 
Under Assumption~\ref{assumption:setting}, suppose that for some $\gamma > 1/2$ and $\nu>1/2$,
\begin{subequations}
\begin{align}
    \label{eq:upper_cs_trunc_U}
    &\|\textnormal{Id}_{\mathcal{U}}-(\mathscr{E}_\infty^{-1}\circ \mathscr{R}_s^*)\circ(\mathscr{R}_s\circ\mathscr{E}_\infty)\|_{L^2_{\varrho}(\mathcal{U};\mathcal{U})}\le c' s^{-\gamma}\quad\text{and}\\
    \label{eq:upper_cs_trunc_V}
    &\|\textnormal{Id}_{\mathcal{V}}-(\mathscr{D}_\infty\circ \mathscr{R}_s^*)\circ(\mathscr{R}_s\circ\mathscr{D}_\infty^{-1})\|_{L^2_{\mathscr{G}_{\sharp}\varrho}(\mathcal{V};\mathcal{V})}\le c' s^{-\nu+1/2}
    \end{align}
\end{subequations}
for all $s\in\mathbb{N}$, with $c'=c'(\mathscr{G},\varrho)>0$ a suitable constant. Further assume that $(\mathscr{E}_{\infty})_\sharp \varrho$ is a quasi-uniform measure over $\mathsf{H}_{0,1}$ and that $(\mathscr{R}_s \circ \mathscr{E}_{\infty})_{\sharp} \varrho$ is absolutely continuous with respect to the uniform measure over $\mathsf{H}_{0,1}$.
Assume that $\mathscr{G}\circ\mathscr{E}_\infty^{-1} \in \mathcal{H}(b)$ for some sequence $b\in\ell^p_{\mathsf{M}}(\mathbb{N})$, where $p\in(0,1)$. For some $\sigma>0$, suppose that $\|e_i\|_{\mathcal{V}}\le\sigma$ for every $i$. Set $d=q=\lceil\tilde{n}\rceil$.
Then with high probability, the minimum in \eqref{eq:training_pb} is attained, and any minimizer $\hat{g}_n\in\mathscr{N}$
yields an operator $\hat{\mathscr{G}}_n\coloneqq \mathscr{D}_{q(n)}\circ \hat{g}_n\circ\mathscr{E}_{d(n)}$ such that
\begin{equation}
\label{eq:upper_cs_bound}
    \|\mathscr{G} - \hat{\mathscr{G}}_n\|_{L^2_{\varrho}(\mathcal{U}; \mathcal{V})} 
    \le c \bigl(
    \tilde{n}^{- \min\{1/p, \gamma, \nu\} + 1/2} + \sigma\bigr)
\end{equation}
for some $c=c(b,\delta,\mathscr{D}_\infty,\mathscr{E}_\infty,\gamma, \nu, p)>0$.
\end{theorem}

Observe that, since $p \in (0,1)$, the first term $\tilde{n}^{- \min\{1/p, \gamma, \nu\}+1/2}$ of the error bound \eqref{eq:upper_cs_bound} is faster than the Monte Carlo rate $n^{-1/2}$ as soon as $\min\{\gamma, \nu\} > 1$. This fast algebraic rate, made possible by the holomorphic regularity of the operator $\mathscr{G}$, is optimal (up to log factors) in a minimax sense; see Thm.~\ref{thm:minimax_holo}. This rate, however, does not immediately imply decay of the generalization error \eqref{eq:upper_cs_bound} unless the noise is absent ($\sigma=0$) or assumed to be vanishing as $n\to\infty$. Notice that the latter scenario is indeed possible: consider, for instance, $\{e_i\}$ arising from numerical PDE solvers whose error can be progressively reduced by considering finer and finer discretizations as the sample size grows.

In line with other bounds proved in the context of function approximation, sometimes referred to as \emph{practical existence theorems} \cite{adcock2024learning, adcock2025near, franco2025practical},
the proof of Thm.~\ref{thm:upper_cs} relies on a construction that emulates a sparse polynomial approximation scheme based on compressed sensing \cite{adcock2022sparse}. This leads to a class of networks $\mathscr{N}$ with a precise and special structure. 
Specifically, while the last layer of a network in $\mathscr{N}$ is a standard fully connected and trainable one, the weights and biases from the first to the second-to-last layer are ``handcrafted''. This means that, during training, they can only take values from a large but finite set. This set corresponds to networks that approximate orthogonal polynomials whose \emph{weighted cardinality} does not exceed a certain threshold \cite[App.~D]{adcock2024optimal}. Notably, the class $\mathscr{N}$ is ``problem agnostic'', i.e., it is independent of the parameter $b$ defining the smoothness of $\mathscr{G}$.

Although the analysis in \cite{adcock2024optimal} considers tanh activations, the same arguments can be readily adapted to other ones, e.g., ReLU. In addition, Thm.~\ref{thm:upper_cs} can be extended to hold for operators between Banach spaces and deliver $L^\infty_\varrho(\mathcal{U}; \mathcal{V})$ error bounds. Finally, similarly to $\tau$ of Thm.~\ref{thm:upper_ep}, the constant $\delta$ of Thm.~\ref{thm:upper_cs} is an artifact of the proof and can be chosen arbitrarily close (but not equal) to zero.

\runinhead{Discussion.}
We conclude this section by illustrating some key differences between Thms.~\ref{thm:upper_ep} and \ref{thm:upper_cs}. First, we stress once more the different nature of these results. Indeed, while Thm.~\ref{thm:upper_ep} relies on a statistical analysis for which the near Monte Carlo rate cannot be expected to be improved in general (see \cite[\S4.1]{reinhardt2024statistical}), Thm.~\ref{thm:upper_cs} is essentially an approximation-theoretic result that allows for faster rates. As such, a key difference is that Thm.~\ref{thm:upper_ep} considers the case of data perturbed by i.i.d.---potentially unbounded---statistical noise, whereas Thm.~\ref{thm:upper_cs} focuses on bounded perturbations. Consequently,
in the first result a nonparametric rate is obtained that is at best of order $n^{-1/2}$. In the second result, a faster-than-Monte Carlo rate is possible when $\gamma$ and $\nu$ in \eqref{eq:upper_cs_trunc_U}--\eqref{eq:upper_cs_trunc_V} are large enough and the noise magnitude satisfies $\sigma \leq c \tilde{n}^{-(1/p-1/2)}$. Moreover, in contrast to Thm.~\ref{th:jakob}, Thm.~\ref{thm:upper_cs} establishes a result in ``high probability''. This means that the probability of the event that the error bound holds is arbitrarily high but fixed (e.g., $0.99$). The dependence of the bounds on such constants is explicitly tracked in \cite{adcock2024optimal}. Still, we note that results in probability can also be obtained in the setting of Thm.~\ref{th:jakob}. Indeed, the analysis in \cite{reinhardt2024statistical}, particularly \cite[Thms.~2.6 and 2.10]{reinhardt2024statistical}, provides concentration inequalities that complement the expectation bounds stated in Thm.~\ref{th:jakob}.

Another key difference concerns the network class $\mathscr{N}$. In Thm.~\ref{thm:upper_ep}, $\mathscr{N}$ contains fully trainable, unconstrained MLPs. In contrast, the class $\mathscr{N}$ in Thm.~\ref{thm:upper_cs} consists of ``handcrafted'' architectures. We note that this gap is partially addressed by \cite[Thm.~3.2]{adcock2024optimal}, which establishes the existence of uncountably many minimizers able to achieve error bounds analogous to \eqref{eq:upper_cs_bound} when $\mathscr{N}$ contains fully trainable MLPs.
Finally, the two classes of holomorphic operators considered in Thms.~\ref{thm:upper_ep} and \ref{thm:upper_cs} differ. In Thm.~\ref{thm:upper_ep}, it is sufficient for $\mathscr{G}\circ\mathscr{E}_{\infty}^{-1}$ to admit a holomorphic extension to an arbitrary open set containing $\mathsf{H}_{r,R}$, whereas in Thm.~\ref{thm:upper_cs}, $\mathscr{G}\circ\mathscr{E}_{\infty}^{-1}$ must be holomorphic on the specific region $\mathcal{R}(b) \supseteq \mathsf{H}_{0,1}$.

\section{Minimax analysis of operator learning}\label{sec:minimax}
This section covers fundamental performance limits of \emph{any} operator reconstruction based on a finite number of data samples. Although there are probabilistic notions of optimal performance \cite{NovakWozniakowski2008}, we only consider worst-case analysis. To this end, let $X$ be a real Banach space of nonlinear operators acting between two Hilbert spaces $\mathcal{U}$ and $\mathcal{V}$.
Define a \emph{method based on $n$ samples} to be any map $\mathsf{M}_n\colon K\subseteq X\to X$ in the set
\begin{align*}
    \mathsf{Map}_n(K;X)\coloneqq \Bigl\{\mathsf{M}\colon f\mapsto \mathsf{D}_n\bigl(f(u_1),\ldots, f(u_n)\bigr) \,\Big|\, \mathsf{D}_n\colon \mathcal{V}^n\to X, \ \{u_j\}_{j=1}^n\subseteq \mathcal{U} \Bigr\},
\end{align*}
where we assume point evaluation to be continuous over $K$ in the norm $\|\cdot\|_{X}$.

The map $\mathsf{M}_n$ factorizes into an abstract \emph{encoding} step, here constrained to be point evaluation of the target operator at $n$ samples, and an abstract \emph{decoding} step, here represented by an arbitrary map from $\mathcal{V}^n$ back to $X$. Note that these encoders and decoders pertain to operators themselves, and not to operator inputs or outputs as in \S\ref{sec:erm}. Still, ERM can be understood as an element of $\mathsf{Map}_n(K;X)$. Based on this setup, define the \emph{minimax reconstruction error based on $n$ samples} to be
\begin{align}\label{eqn:defn_minimax_sampling_width}
    s_n(K)_{X} \coloneqq \inf_{\mathsf{M}\in \mathsf{Map}_n(K;X)}\sup_{f\in K}\|f-\mathsf{M}(f)\|_{X}.
\end{align}
Eqn.~\eqref{eqn:defn_minimax_sampling_width} is also called the \emph{nonlinear sampling $n$-width} of $K$ in $X$ \cite{kovachki2024data}. Other notions of $n$-width are also possible, such as those based on general encoding \cite{kovachki2024data}, linear decoding only \cite{siegel2024sharp}, or random samples only \cite{parhi2025upper}.
The number $s_n(K)_X$ records the worst error over the set of targets in $K$ incurred by the best method in $\mathsf{Map}_n(K;X)$. The present paper focuses on the \emph{rate of convergence} of $ s_n(K)_{X}\to 0$ as $n\to\infty$ for particular choices of the target operator class $K$ and ambient space $X$. For $\varrho$ with finite second moment, we fix $X\coloneqq L^2_\varrho(\mathcal{U};\mathcal{V})$. It is convenient to let $\{\lambda_j(\mathrm{Cov}(\varrho))\}_{j\in\mathbb{N}}$ denote any nonincreasing rearrangement of the eigenvalues of the covariance operator $\mathrm{Cov}(\varrho)$ of $\varrho$. In the remainder of this section, we assume that $\mathrm{Cov}(\varrho)$ is strictly positive definite.

\runinhead{Lipschitz and Frech\'et differentiable operators.}
We begin with a negative result about the sample efficiency of operator learning.
For $k\in\mathbb{N}$, let $C^k(\mathcal{U};\mathcal{V})$ denote the space of $k$-times Frech\'et differentiable operators with all derivatives uniformly bounded over $\mathcal{U}$ \cite[App.~A]{kovachki2024data}. The following hardness result follows from \cite[Thm.~2.12]{kovachki2024data}.\footnote{Perusal of the proof of \cite[Thm.~2.12]{kovachki2024data} shows that \cite[Thm.~2.12, Eqn.~(2.4)]{kovachki2024data} is actually valid for covariance eigenvalues satisfying $\lambda_j\gtrsim j^{-2\alpha}$, not just $\lambda_j\gtrsim j^{-\alpha}$ as claimed in \cite[\S2.3.2]{kovachki2024data}.}
\begin{theorem}\label{thm:lower_Ck}
    Let $\varrho$ be a Gaussian measure on $\mathcal{U}$. Suppose that for some $\vartheta>1/2$, it holds that $\inf_{j\in\mathbb{N}}j^{2\vartheta}\lambda_j(\mathrm{Cov}(\varrho))>0$. Fix $k\in \mathbb{N}$. If $\mathcal{V}$ is a Hilbert function space containing constant functions, then there exists $c=c(k,\vartheta)>0$ such that for all $n\in\mathbb{N}$,
    \begin{align*}
        s_n\bigl(\bigl\{\mathscr{G}\in C^k(\mathcal{U};\mathcal{V})\,\big|\, \|\mathscr{G}\|_{C^k}\leq 1 \bigr\}\bigr)_{L^2_\varrho(\mathcal{U};\mathcal{V})}\geq c\bigl(1+\log(n)\bigr)^{-k(\vartheta+3)}.
    \end{align*}
\end{theorem}
Thm.~\ref{thm:lower_Ck} shows that, for $K$ equal to the $C^k$ unit ball and no matter how large the order $k$ of Frech\'et differentiability is, the nonlinear sampling $n$-width of $K$ in $L^2_{\varrho}$ decays no faster than polylogarithmically; cf.~\cite{cheng2026learning,liu2024deep,liu2024neural} for upper bounds. Therefore, \emph{no method} can exhibit algebraic sample complexity uniformly over $K$. This is a ``curse of sample complexity'' for operator learning.
Although \cite[Thm.~2.12]{kovachki2024data} is stated only for $C^k(\mathcal{U};\mathbb{R})$ functionals, we can identify such a functional with a $C^k(\mathcal{U};\mathcal{V})$ operator taking value as a constant function. This fact implies that the sampling $n$-width of $C^k$ functionals is bounded above by the sampling $n$-width of $C^k$ operators, leading to Thm.~\ref{thm:lower_Ck}.
Moreover, Thm.~\ref{thm:lower_Ck} with $k=1$ implies a sampling $n$-width lower bound for $K$ equal to the unit Lipschitz norm ball because this set contains a $C^1$ norm ball.
Related work refines Thm.~\ref{thm:lower_Ck} in the $k=1$ case to allow for unbounded Lipschitz operators, general eigenvalue decay of $\mathrm{Cov}(\varrho)$, and general linear measurements \cite[\S5]{adcock2024lip}; the conclusions are qualitatively the same.

\runinhead{Holomorphic operators.}
We now show how assuming stronger regularity of the operator allows for faster minimax rates. Let $\iota\colon \mathcal{U} \to \mathbb{R}^{\mathbb{N}}$ measurable be such that $\iota_\sharp \varrho$ is a quasi-uniform measure supported on $\mathsf{H}_{0,1} = [-1,1]^{\mathbb{N}}$. Suppose that $\iota|_{\mathrm{supp}(\varrho)}\colon \mathcal{U} \to \ell^{\infty}(\mathbb{N})$ is Lipschitz. Consider the holomorphic class $\mathcal{H}(b) \circ \iota \coloneqq \{f\circ \iota \,|\, f \in \mathcal{H}(b)\}$, where $\mathcal{H}(b)$ is as in \eqref{eq:def_H(b)}. Then \cite[Thm.~4.1]{adcock2024optimal} shows the following.
\begin{theorem}
\label{thm:minimax_holo}
For any $p \in (0,1)$, there exists a constant $c(p)>0$ such that the following hold.  First, for each $n \in \mathbb{N}$, there exists $b = b^{(n)}\in \ell^p_{\mathsf{M}}(\mathbb{N})$, with $b^{(n)}_j\geq 0$ for each $j$ and $\|b^{(n)}\|_{p,\mathsf{M}}=1$, such that $s_n(\mathcal{H}(b^{(n)})\circ \iota)_{L^2_{\varrho}(\mathcal{U};\mathcal{V})} \geq c(p)\, n^{-(1/p-1/2)}$.
Second, there exists $b\in \ell^p_{\mathsf{M}}(\mathbb{N})$, with $b_j\geq 0$ for each $j$ and $\|b\|_{p,\mathsf{M}}=1$, such that
\begin{align*}
    s_n\bigl(\mathcal{H}(b)\circ \iota\bigr)_{L^2_{\varrho}(\mathcal{U};\mathcal{V})} \geq c(p) \, n^{-\left(\frac{1}{p}-\frac{1}{2}\right)}\log^{-\frac{2}{p}}(2n) \quad \text{for all}\quad n \in \mathbb{N}.
\end{align*}
\end{theorem}
Thm.~\ref{thm:minimax_holo} shows that the rate $n^{-(1/p-1/2)}$ is optimal, up to log factors, for the class of holomorphic operators $\mathcal{H}(b)\circ \iota$. Notably, this rate is achieved by the empirical risk minimizer from Thm.~\ref{thm:upper_cs} when $\min\{\gamma, \nu\} \geq 1/p$ and $\sigma =0$.

\runinhead{Neural-based classes of target operators.}
In view of the hardness of operator learning over the $k$-times Frech\'et differentiable class from Thm.~\ref{thm:lower_Ck}, it is natural to ask whether it is possible to achieve better sample complexity over a smaller class $K\subset C^k(\mathcal{U};\mathcal{V})$ of operators. We have seen in Thm.~\ref{thm:minimax_holo} that taking $K\coloneqq \mathcal{H}(b) \circ \iota$ suffices to achieve arbitrarily fast algebraic rates. However, the extreme smoothness that holomorphy imparts is a strong condition. A different line of work defines a class $K^\alpha$ of target operators that are efficiently approximated at rate $\alpha$ by specific operator learning architectures \cite{grohs2025theory,kovachki2024data}, such as the DeepONet \cite{lu2021learning} or FNO \cite{li2020fourier}. Although there are results for DeepONets \cite[\S3.4]{grohs2025theory}, we focus on the FNO as a case study.

Fix $\Omega\coloneqq (0,1)^{d_{\Omega}}$ and $\mathcal{U}=\mathcal{V}\coloneqq L^2(\Omega)$. To define $K^\alpha\coloneqq K^\alpha_{\mathsf{FNO}}$, we first recall that the FNO is a kernel integral neural operator architecture formed as a composition of $L$ layers \cite{boulle2024mathematical,kovachki2024operator,li2020fourier}. These layers include integral operators whose convolutional kernels are parametrized by a finite number of Fourier coefficients indexed by wavenumbers $k\in\mathbb{Z}^{d_\Omega}$ satisfying $\|k\|_{\ell^\infty}<\kappa$ for some cutoff $\kappa\in \mathbb{N}$. Each hidden layer maps between latent function spaces with codomain dimension $d_{\mathrm{c}}$. We use ReLU activations as in \cite{grohs2025theory,kovachki2024data}. For any $m\in\mathbb{N}$, define $\mathscr{N}_m^{\mathsf{FNO}}$ to be the set of all FNOs $\Psi=\Psi_\theta\colon\mathcal{U}\to\mathcal{V}$ with tunable parameters $\theta$ such that $\max\{\kappa^{d_\Omega}, d_{\mathrm{c}}, L\}\leq m$ and $\|\theta\|_{\ell^\infty}\leq \exp(m)$. Next, for any compact set $\mathcal{K}\subset\mathcal{U}$, let $e_\mathcal{K}(f, \mathscr{N})\coloneqq \inf_{\Psi\in\mathscr{N}}\|f-\Psi\|_{C(\mathcal{K};\mathcal{V})}$ and 
\begin{align*}
    \Gamma^\alpha_{\mathsf{FNO}}(f)_{\mathcal{K}}\coloneqq \max\Bigl(\|f\|_{C(\mathcal{K};\mathcal{V})}, \, \sup_{m\in\mathbb{N}}m^\alpha e_\mathcal{K}(f, \mathscr{N}_m^{\mathsf{FNO}})\Bigr).
\end{align*}
Finally, define the \emph{space of target operators well approximated by FNOs} to be
\begin{align}\label{eqn:fno_space}
    K^\alpha_{\mathsf{FNO}}(\mathcal{K})\coloneqq \bigl\{\mathscr{G}\in C(\mathcal{K};\mathcal{V})\,\big|\, \inf\{t>0\,|\, \Gamma^\alpha_{\mathsf{FNO}}(\mathscr{G}/t)_{\mathcal{K}}\leq 1\}\leq 1 \bigr\}.
\end{align}
To show that $K^\alpha$ in \eqref{eqn:fno_space} escapes the curse of sample complexity, we study upper bounds on $s_n(K^\alpha)_{L^2_\varrho}$. To this end, define the \emph{optimal algebraic minimax rate exponent} to be
\begin{align*}
    \beta^\star(K)_X\coloneqq \sup\Bigl\{\beta\geq 0 \, \Big|\, \sup_{n\in\mathbb{N}}n^\beta s_n(K)_X<\infty\Bigr\}.
\end{align*}
Notice that $K_1\subseteq K_2$ implies $\beta^\star(K_1)_X\geq \beta^\star(K_2)_X$.
We remark that \cite[Thm.~1]{parhi2025upper}, combined with \cite[Thm.~10]{siegel2024sharp}, gives a general (albeit loose) recipe to lower bound $\beta^\star(K)_X$ based on the $L^\infty$ metric entropy of $K$.

From \cite[Thm.~3.3]{kovachki2024data} and \cite[Thm.~3.21]{grohs2025theory}, we deduce the next theorem.
\begin{theorem}\label{thm:fno_optimal_rate}
Let $\Omega\coloneqq (0,1)^{d_{\Omega}}$. Let $\varrho$ be supported on a compact set $\mathcal{K}\subset L^2(\Omega)$. Suppose that $\mathcal{K}$ is convex and is not contained in any finite-dimensional subspace of $L^2(\Omega)$. Then for any $\alpha>0$, it holds that
    \begin{align}\label{eqn:fno_optimal_rate}
        \frac{1}{2}\biggl(\frac{1}{1+8/\alpha}\biggr)\leq  \beta^\star\bigl(K_{\mathsf{FNO}}^{\alpha}(\mathcal{K})\bigr)_{L^2_\varrho\left(L^2(\Omega);L^2(\Omega)\right)}\leq \frac{1}{2}.
    \end{align}
\end{theorem}
The lower bound on the optimal exponent follows from an algebraic upper bound on $s_n(K^\alpha_{\mathsf{FNO}})_{L^2_\varrho}$ based on ERM \cite[Thm.~3.3]{kovachki2024data}. For the upper bound on the optimal exponent, \cite[Rem.~3.19]{grohs2025theory} (taking $\kappa=1$) and perusal of the proof of \cite[Lem.~3.22]{grohs2025theory} show that \cite[Eqn.~(3.30)]{grohs2025theory} remains valid with $K_{\mathsf{FNO}}^{\alpha}(\mathcal{K})$ in place of $\mathsf{U}_{\ell,\mathsf{NO}}^{\alpha,\infty}$ in that lemma. Then \cite[Eqn.~2.6]{grohs2025theory} and the proof of \cite[Thm.~3.21]{grohs2025theory} yields \eqref{eqn:fno_optimal_rate}.

The hypotheses on $\varrho$ in Thm.~\ref{thm:fno_optimal_rate} are satisfied by invoking \cite[Prop.~3.5]{grohs2025theory}; see \cite[Eqn.~(4.22), pp.~36--37]{de2025extension} for an explicit example.
Thm.~\ref{thm:fno_optimal_rate} shows that the optimal algebraic minimax rate exponent is asymptotically $1/2$ as $\alpha\to\infty$. 
However, for fixed $ \alpha<\infty $, the sharp value of the exponent remains open. The lower estimate in \eqref{eqn:fno_optimal_rate} is due to an ERM construction that may be suboptimal. Alternatively, the deterioration as $ \alpha $ decreases may be intrinsic because $ K_{\mathsf{FNO}}^{\alpha}(\mathcal{K}) $ becomes larger for smaller $ \alpha $.

\runinhead{Statistical noise models.}
Thus far, this section has worked in the ``optimal recovery'' setting in which the data consist of exact point values of the target operator. To address observations corrupted by random noise, e.g., as in Thm.~\ref{th:jakob}, let $\eta$ be a centered probability measure on $\mathcal{V}$ and define the \emph{statistical nonlinear sampling $n$-width} $\tilde{s}_n^\eta(K;\sigma)_{X}$ by
\begin{align}\label{eqn:noisy_width}
    \inf_{(\{u_j\}_{j=1}^n, \mathsf{D}_n)}\sup_{f\in K}\mathbb{E}_{\{\xi_j\}_{j=1}^n\sim\eta^{\otimes n}}\bigl\|f-\mathsf{D}_n\bigl(f(u_1)+\sigma\xi_1,\ldots, f(u_n)+\sigma\xi_n\bigr)\bigr\|_{X}.
\end{align}
Here $\sigma>0$ is the noise level. Although in \eqref{eqn:noisy_width} the design points in the infimum are fixed, it is also common to consider random design and minimize over the distribution of the random points \cite{adcock2025towards}.
Writing $L^2_\varrho\coloneqq L^2_\varrho(\mathcal{U};\mathcal{V})$, \cite{adcock2025towards} considers $\tilde{s}_n^\eta(K;\sigma)_{L^2_\varrho}$ with $K\coloneqq K_{\mathsf{Lip}}$ equal to the unit Lipschitz norm ball and possibly non-Gaussian $\varrho$.
A sharp characterization of $\tilde{s}_n^\eta(K;\sigma)_{L^2_\varrho}$ in terms of $n$ is given when the eigenvalues of $\mathrm{Cov}(\varrho)$ decay exponentially \cite[Thm.~2.5]{adcock2025towards}. For algebraic decay, a similar result implies that---up to the exponent---the inverse polylogarithmic behavior of Thm.~\ref{thm:lower_Ck} is essentially sharp.
\begin{theorem}\label{thm:lip_upper_lower}
    Let $K_{\mathsf{Lip}}$ denote the Lipschitz norm unit ball. Instate the hypotheses of Thm.~\ref{thm:lower_Ck}. If $\lambda_j(\mathrm{Cov}(\varrho))\asymp j^{-2\vartheta}$ and $\eta$ is Gaussian, then for all sufficiently large $n$,
    \begin{align}\label{eqn:lip_upper_lower}
        \bigl(\log(n)\bigr)^{-(\vartheta+3)}\lesssim s_n(K_{\mathsf{Lip}})_{L^2_\varrho}\leq \tilde{s}_n^\eta(K_{\mathsf{Lip}};\sigma)_{L^2_\varrho}\lesssim \biggl(\frac{\log\bigl(n\sigma^{-2}\bigr)}{\log\log\bigl(n\sigma^{-2}\bigr)}\biggr)^{-(\vartheta-1/2)}.
    \end{align}
\end{theorem}
The rightmost inequality in \eqref{eqn:lip_upper_lower} is from \cite[Thm.~2.9]{adcock2025towards}. The first follows from Thm.~\ref{thm:lower_Ck}. The second inequality is due to the fact that the noise-free sampling $n$-width \eqref{eqn:defn_minimax_sampling_width} is bounded above by the noisy one. Indeed, for any $\mathsf{D}_n\colon \mathcal{V}^n\to X$, define $\tilde{\mathsf{D}}_n(v_1,\ldots, v_n)\coloneqq\mathbb{E}_{\{\xi_j\}_{j=1}^n\sim\eta^{\otimes n}}[\mathsf{D}_n(v_1+\sigma\xi_1,\ldots, v_n+\sigma\xi_n)]$. Then $\tilde{\mathsf{D}}_n$ and any $\{u_j\}_{j=1}^n\subset \mathcal{U}$ define an element of $\mathsf{Map}_n(K;X)$ that is feasible for \eqref{eqn:defn_minimax_sampling_width}. Consequently, Jensen's inequality for the $X$ norm and the definition \eqref{eqn:noisy_width} imply that $s_n(K)_{X}\leq \tilde{s}_n^\eta(K;\sigma)_{X}$.

\runinhead{Discussion.}
This section shows that algebraic sample complexity for $L^2_\varrho$-learning of Lipschitz (or even $C^k$) operators is impossible in the worst case. On the other hand, minimax errors for holomorphic operators can achieve any desired algebraic rate if the summability exponent $p$ is sufficiently close to zero. As a middle ground, the sample complexity over $K^\alpha_{\mathsf{FNO}}$ is algebraic, but can be no better than $n^{-1/2}$ even if $\alpha\gg 1$. Taken together, the results in \S\ref{sec:erm} and \S\ref{sec:minimax} imply that
$ \max\{s_n(\mathcal{H}(b)\circ \iota)_{L^2_{\varrho}},\, s_n(K^\alpha_{\mathsf{FNO}})_{L^2_\varrho}\}\lesssim s_n(K_{\mathsf{Lip}})_{L^2_\varrho}$.
A further separation of the left-hand side would require lower bounds on $s_n(K^\alpha_{\mathsf{FNO}})_{L^2_\varrho}$ or showing that elements of $\mathcal{H}(b)\circ \iota$ are well approximated by FNOs.

\section{Conclusions and open problems}\label{sec:conclusion}
This paper reviews recent progress in the theory of operator learning by presenting two convergence analyses of ERM for holomorphic operators (Thms.~\ref{thm:upper_ep}~and~\ref{thm:upper_cs} in \S\ref{sec:erm}) and by surveying multiple results about minimax sample complexity rates (Thms.~\ref{thm:lower_Ck}--\ref{thm:lip_upper_lower} in \S\ref{sec:minimax}). Several open problems naturally arise from this paper.

Based on the results presented in \S\ref{sec:erm}, it is natural to ask whether faster-than-Monte Carlo rates are possible when considering ERM with \emph{fully trainable networks}, especially in the absence of noise. Currently, an answer is elusive because the arguments used to prove Thms.~\ref{thm:upper_ep}~and~\ref{thm:upper_cs} do not easily extend to such a setting. In fact, the technical approach underlying Thm.~\ref{thm:upper_ep} is likely to be suboptimal when $\sigma=0$; conversely, the construction considered in Thm.~\ref{thm:upper_cs} is still not general enough to accommodate fully trainable architectures. Notably, these issues are not intrinsic to operator learning, but rather constitute an open question even in finite dimensions. A remaining question is whether the nonparametric rate in Thm.~\ref{thm:upper_ep} (which is at best $n^{-1/2}$) is an artifact of the statistical learning analysis in the noisy setting, or a genuine minimax barrier.

One interesting future research direction in minimax analysis concerns the random noise setting from \S\ref{sec:minimax}. For classical function approximation, it is known that $s_n(K)_{X}$ from \eqref{eqn:defn_minimax_sampling_width} and $\tilde{s}_n^\eta(K;\sigma)_{X}$ from \eqref{eqn:noisy_width} behave quite differently as a function of $n$, although they do agree in the limit $\sigma\to 0$ \cite{devore2025optimal}. An open problem is to obtain a similar \emph{noise level aware} characterization of minimax rates for operator learning. Another important challenge is to identify classes of operators that are practically relevant in scientific applications, but still enjoy \emph{algebraic} minimax sample complexity.

Collectively, these considerations open up several promising research directions at the interface of the topics in \S\ref{sec:erm} and \S\ref{sec:minimax}.
First, can one prove sharp bounds on the statistical $n$-width $\tilde{s}_n^\eta(K;\sigma)_{L^2_\varrho}$ \eqref{eqn:noisy_width} for $K$ the holomorphic classes considered in Thms.~\ref{thm:upper_ep}--\ref{thm:upper_cs}? This will help clarify whether faster-than-Monte Carlo rates are possible or not under suitable statistical noise models. Moreover, is ERM optimal in the minimax sense of \eqref{eqn:defn_minimax_sampling_width} for target classes beyond holomorphic operators? Other sets of interest could include $C^\infty$ operators or those arising from specific parametric PDE problems. In addition, as suggested by the results in \S\ref{sec:minimax}, an interesting alternative could be to consider architecture-oriented operator classes, in contrast to those based on more classical regularity assumptions.

\begin{acknowledgement}
S.B.\ acknowledges support from NSERC through grant RGPIN-2020-06766
and FRQ - Nature et Technologies through grant 359708. N.R.F.\ acknowledges support from project DREAM (Reduced Order Modeling and Deep Learning for the real-time approximation of PDEs), grant no. FIS00003154, funded by the Italian Science Fund (FIS) and by Ministero dell'Università e della Ricerca (MUR); his present research is part of the activities of project Dipartimento di Eccellenza 2023-2027, Department of Mathematics, Politecnico di Milano, funded by MUR.
\,N.H.N.\ acknowledges support from a Klarman Fellowship through Cornell University's College of Arts \& Sciences.
The authors are also grateful to Gregor Maier and an anonymous referee for helpful comments.
\end{acknowledgement}

\bibliographystyle{abbrv}
\bibliography{references}

\end{document}